\documentclass{article}

\usepackage{latexsym}
\usepackage{a4wide}
\usepackage{amscd}
\usepackage{graphics}
\usepackage{amsmath}
\usepackage{amssymb}
\usepackage{mathrsfs}
\input xy
\xyoption{all}

\newcommand{\N}{\mathbb{N}}                     
\newcommand{\Z}{\mathbb{Z}}                     
\newcommand{\R}{\mathbb{R}}                     
\newcommand{\set}[2]{\left\{{#1}\mid{#2}\right\}}       
\newcommand{\qed}{\hfill $\Box$ \bigskip}       
\newcommand{\im}{\mathrm{Im\,}}                 
\newcommand{\dist}{\mathrm{dist\,}}             
\newcommand{\ind}{\mathrm{ind\,}}               
\newcommand{\codim}{\mathrm{codim}}           
\newcommand{\rank}{\mathrm{rank\,}}             
\newcommand{\crit}{\mathrm{crit}}               
\newcommand{\rest}{\mathrm{rest}\,}             

\newtheorem{thm}{\sc Theorem}[section]      
\newtheorem{lem}[thm]{\sc Lemma}            
\newtheorem{rem}[thm]{\sc Remark}           

\title{A Morse complex for Lorentzian geodesics}

\author{Alberto Abbondandolo and Pietro Majer}

\date{Universit\`a di Pisa}
 
\begin{document}

\maketitle

\begin{abstract}
We prove the Morse relations for the set of all geodesics connecting two
non-conjugate points on a class of globally hyperbolic Lorentzian
manifolds. We overcome the difficulties coming from the fact that the 
Morse index of every geodesic is infinite, and from the lack of the 
Palais-Smale condition, by using the Morse complex approach.  
\end{abstract}

\renewcommand{\theenumi}{\roman{enumi}}
\renewcommand{\labelenumi}{(\theenumi)}

\section*{Introduction}

Let $M$ be a smooth connected manifold without boundary of dimension 
$n+1$, and let $h$ be a Lorentzian structure on $M$: this means that $h$ is a
non-degenerate symmetric (0,2)-tensor on $M$ having $n$ positive
eigenvalues and one negative eigenvalue (see \cite{one83} and
\cite{bee96} for 
foundational results on Lorentzian geometry). The Lorentzian structure 
$h$ induces a unique Levi-Civita covariant derivative $\nabla$ on $M$, 
and a geodesic on $(M,h)$ is a curve $\gamma: I \rightarrow M$ whose
velocity $\gamma'$ is parallelely transported with respect to
this covariant derivative. 

We are interested in the problem of classifying all geodesics
connecting two fixed points $z_0,z_1$ in $M$. These geodesics are
critical points of the energy functional
\[
E(\gamma) = \frac{1}{2} \int_0^1 h(\gamma'(t),
\gamma'(t)) \, dt,
\]
defined on the Hilbert manifold $\mathscr{M}$ consisting of the 
curves $\gamma:[0,1] \rightarrow M$ connecting $z_0$ to $z_1$, of
Sobolev class $W^{1,2}$. 

In the Riemannian case, the geodesics $\gamma$ connecting $z_0$ to
$z_1$ can be classified according to their Morse index $i(\gamma)$,  
that is the number of negative eigenvalues of the second differential 
of the Riemannian energy functional at $\gamma$. 
This index has also other geometrical
interpretations: it coincides with $i_\mathrm{con}(\gamma)$, 
the number of conjugate points along $\gamma$, counted with
multiplicity, and with the Maslov index $i_{\mathrm{Maslov}}(\gamma)$
of a path of Lagrangian subspaces obtained by passing to the
Hamiltonian formulation and linearizing the geodesic flow along
$\gamma$. When the Riemannian manifold $M$ is complete and when
$z_1$ is not conjugated to $z_0$ along any geodesic, the geodesics
connecting $z_0$ to $z_1$ satisfy the following Morse
relations. Denote by $c_k\in \N \cup \{+\infty\}$ the number of
geodesics connecting $z_0$ to $z_1$ with Morse index $k$. Then there
exists a formal series $Q$ with coefficients in $\N \cup \{+\infty\}$
such that
\begin{equation}
\label{morserel}
\sum_{k=0}^{\infty} c_k x^k = \sum_{k=0}^{\infty} \rank H_k(\Omega(M))
\, x^k + (1+x) Q(x),
\end{equation}
where $\Omega(M)$ is the space of based loops on $M$, and $H_k$
denotes singular homology with integer coefficients (for a modern
presentation of this famous result by Morse see \cite{pal63}). 
Equivalently, the Morse relations (\ref{morserel}) can be restated in
the following way:
denoting by $C_k$ the free Abelian group generated by the
geodesic connecting $z_0$ to $z_1$ with Morse index $k$, there exists
homomorphisms $\partial_k: C_k \rightarrow C_{k-1}$ such that
$\{C_k,\partial_k\}_{k\in \N}$ is a chain complex with homology isomorphic to
the singular homology of $\Omega(M)$. Actually, the chain complex 
$\{C_k,\partial_k\}_{k\in \N}$ can be defined as the cellular chain
complex associated to a cellular filtration of $\mathscr{M}$ produced
by the negative gradient flow of the Riemannian energy functional (see for instance \cite{kli78} or
\cite{ama06m}).  

In the Lorentzian case, all the critical points of the energy
functional have infinite Morse index. However, a finite relative index
$i_{\mathrm{rel}}(\gamma)$ can be defined in the following way. We fix 
a distribution $W$ of hyperplanes in $TM$ on which $h$ is positive
definite, and we define $\mathscr{W}$ as the subbundle of $T\mathscr{M}$
such that $\mathscr{W}(\gamma)$ is the space of $W^{1,2}$ sections of
$\gamma^*(TM)$ with values in $W$. Then the relative Morse index of a
critical point $\gamma$ of $E$ on $\mathscr{M}$ is the integer
\[
i_{\mathrm{rel}} (\gamma) = \ind (V^-(\nabla^2
E(\gamma)),\mathscr{W}(\gamma)),
\]
where $V^-(\nabla^2 E(\gamma))$ is the negative eigenspace of the
Hessian of $E$ at $\gamma$, and $\ind(\mathscr{V},\mathscr{W})$
denotes the Fredholm index of the pair $(\mathscr{V},\mathscr{W})$,
that is the integer $\dim \mathscr{V} \cap \mathscr{W} - \codim
(\mathscr{V} + \mathscr{W})$.   

On the other hand, the Maslov index $i_{\mathrm{Maslov}} (\gamma)$ of the 
geodesic $\gamma$ is still a well defined integer, which can now be
also negative. Musso, Pejsachowicz, and Portaluri \cite{mpp05} have shown that
an index $i_{\mathrm{con}} (\gamma)$, counting the
conjugate points along $\gamma$, can still be defined: the difficulty 
given by the fact that in the
case of a non-positive $h$ conjugate points may accumulate can be 
overcome by defining $i_{\mathrm{con}}$
as the topological degree of a map obtained by extending the equation
for Jacobi fields to the complex plane. Results by
these authors together with previous results by Piccione and Tausk 
\cite{pt02}, imply that these indices coincide\footnote{More
  generally, this is true for an arbitrary semi-Riemannian manifold.}:
\begin{equation}
\label{ri=i=i}
i_{\mathrm{rel}}(\gamma) = i_{\mathrm{con}}(\gamma) = i_{\mathrm{Maslov}}
(\gamma).
\end{equation}

As in the Riemannian case, this relative Morse index plays an
important role in the local bifurcation theory of geodesics, see
for instance \cite{ppt04} and \cite{pp05}. 

A global Morse theory for 
{\em time-like} geodesics (that is geodesics $\gamma$ such that
$h(\gamma',\gamma')<0$) was developed by Uhlenbeck \cite{uhl75} for a
class of {\em time oriented globally hyperbolic} Lorentzian manifolds:
time orientability means that there exists a smooth vector field $v$ on $M$ 
such that $h(v,v)<0$, while global hyperbolicity requires that for
every pair of points $z_0,z_1$ in $M$ the space of time-like geodesics
joining them is pre-compact in the compact-open topology. Geroch has
discovered in \cite{ger70} that a time oriented 
globally hyperbolic manifold $M$ is isometric to a product 
$X \times \R$, $X$ an $n$-dimensional
manifold, with Lorentzian structure of the form
\begin{equation}
\label{hyperbolic}
h = \alpha_{ij}(x,y)dx^i \otimes dx^j - \beta(x,y)dy^2,
\end{equation}
where $(x,y)=(x^1,\dots,x^n,y)$ is a local coordinate system in $X \times \R$, 
the matrix $\alpha_{ij}(x,y)$ is positive definite and $\beta$ is a
positive function (see \cite{bs03} and \cite{bs05} for a rigorous proof). 
Then Uhlenbeck's theory makes use of the fact that
time-like geodesics locally maximize the Lorentzian length, so that
methods from standard Morse theory can be applied. Uhlenbeck's results imply that the time-like geodesics joining two given non-conjugate points $z_0=(x_0,y_0)$ and $z_1=(x_1,y_1)$ satisfy Morse relations analogous to (\ref{morserel}), where the homology of the based loop space is replaced by the homology of the space of all {\em time-like curves} joining $z_0$ and $z_1$. The homology of this space depends on the end-points $z_0=(x_0,y_0)$ and $z_1=(x_1,y_1)$, and as  $y_1$ tends to $+\infty$ it converges to the homology of the based loop space $\Omega(X)$, so the Riemannian Morse relations are somehow recovered in the limit of two points infinitely far in time. 

The aim of his paper is to show that if one considers {\em all}  geodesics joining $z_0$ and $z_1$, full Morse relations hold without taking the limit for $y_1\rightarrow +\infty$.  Here is the main result of this paper:

\bigskip 

\noindent {\sc Theorem. }
{\em Let $M=X \times \R$ be a connected manifold endowed with the
Lorentzian structure $h$ defined by (\ref{hyperbolic}). 
Assume that $X$ is compact and that there is some $s_0$ such that
$X \times ]-s,s[$ is convex for every $s>s_0$. 
Let $z_0$ and $z_1$ be two points in $X \times
]-s_0,s_0[$, assumed to be non-conjugate along every geodesic joining
them. Given $k\in \Z$, let $C_k$ be the free Abelian group generated 
by the geodesics $\gamma$ joining $z_0$ and $z_1$ of relative index $k$. Then
there are boundary homomorphisms $\partial_k : C_k \rightarrow C_{k-1}$ such
the homology of the chain complex $\{C_k,\partial_k\}_{k\in \Z}$ is isomorphic
to the singular homology of $\Omega(X)$, the based loop space of $X$.}

\bigskip

We recall that
an open subset $\Omega\subset M$ is said to be {\em convex} if every
geodesic $\gamma:]-\epsilon,\epsilon[\rightarrow M$ with $\gamma(0)\in
\partial \Omega$ and $\gamma(]-\epsilon,\epsilon[)\subset
\overline{\Omega}$ satisfies $\gamma(]-\epsilon,\epsilon[)\subset
\partial \Omega$. The convexity assumption required in the above
result is equivalent to asking
the matrix $(\partial_y \alpha_{ij}(x,y))$ to be non-positive for
$s>s_0$ and non-negative for $s<-s_0$. Once $z_0$ is fixed, the assumption
that $z_1$ should be non-conjugate  to $z_0$ along any geodesic holds
for a residual set of points $z_1$ (see \cite{uhl75}, Theorem 1 (a)).
The above theorem holds also if we replace the
compactness hypothesis on $X$ by suitable bounds on the coefficients of
the Lorentzian structure (see Theorem \ref{main} below). Note that, although hyperbolic Lorentzian manifolds may have geodesics with negative relative index, the contribution to homology of such geodesics is null.

A Morse theory of this kind involving all geodesics
connecting two non-conjugate points has been developed by Masiello (see
\cite{mas94}, Chapter 4) for a class of {\em standard stationary} 
Lorentzian manifolds, that is in the case of a Lorentzian structure 
of the form
\[
h = \alpha_{ij}(x) dx^i \otimes dx^j + \delta_i(x) (dx^i \otimes dy +
dy \otimes dx^i) - \beta(x) dy^2.
\]   
This is a much smaller class than the class of time oriented globally 
hyperbolic manifolds. Given some curve $x:I \rightarrow X$, the
equation for the $y$-component of a geodesic can be uniquely solved,
and solving the geodesic equation is equivalent to solving a second
order ODE for $x$ associated to a nice coercive functional. As a
consequence, on a standard stationary Lorentzian manifold all the
geodesics have non-negative index, and the homology of $\Omega(X)$ is recovered.  

Arbitrary geodesics on globally hyperbolic Lorentzian manifolds of the kind considered in the above theorem were studied by Benci, Fortunato, and Masiello in \cite{bfm94} (see also \cite{mas94}, Chapter 8). By using finite dimensional reductions on
the $y$-component of the geodesic equation, they 
proved the existence of at least one geodesic connecting $z_0$ to $z_1$. Actually, their proof could be refined to prove the weak Morse relations, that is the inequalities $c_k \geq \rank H_k(\Omega(X))$ for every $k\in \N$. See also the review paper \cite{abfm03}. 

The strong Morse relations stated in our main theorem seem to go
beyond the methods used in \cite{bfm94}. One of the main difficulties 
in this problem is that the energy
functional of a globally hyperbolic Lorentzian manifold may fail to satisfy the
Palais-Smale condition at every energy level (unlike the stationary
case). This problem is overcome
in \cite{bfm94} by proving an a priori estimate on the $C^1$ norm of 
the second component of any curve $\gamma=(x,y)$ solving $\nabla
E(x,y)=\lambda y$ for some $\lambda \geq 0$, with energy
$E(\gamma)$ bounded above. Here we need to sharpen this estimate, by 
proving that
\begin{equation}
\label{apriori}
\nabla E(x,y)=\lambda y \mbox{ with } \lambda\geq 0 \quad \Rightarrow \quad
\|y\|_{C^1}^2 \leq p \max\{E(x,y),0\} + q,
\end{equation}
for suitable constants $p,q$ (see Lemma
\ref{stima} below)\footnote{Actually, 
in \cite{bfm94} the
  convexity assumption is replaced by a sort of convexity at
  infinity. This is a weaker assumption which does not imply
  (\ref{apriori})}. Here is where the convexity assumption is used.

Then our argument uses in an essential way the Morse complex
approach developed
in \cite{ama05,ama08p}, instead than classical Morse theory applied to 
finite dimensional reductions.

The idea is to define the boundary homomorphism $\partial_k:C_k
\rightarrow C_{k-1}$ by counting the negative gradient flow lines of
$E$ connecting critical points of relative Morse index difference
one\footnote{This is the basic idea of Floer homology, but in a 
different contest. Here indeed we are dealing with a well-defined
gradient flow on an infinite dimensional Hilbert manifold, whereas
Floer homology deals with an elliptic PDE (see \cite{flo89a}).}. 
The a priori estimate (\ref{apriori}) allows to build a 
pseudo-gradient vector field for $E$ such that all the flow lines
connecting two critical points are bounded. Since the Palais-Smale
condition holds on bounded subsets of $\mathscr{M}$, one gets enough
compactness to define the boundary operator $\partial_k$. Here weaker
a priori estimates as the one proved in \cite{bfm94} would be sufficient, but
the estimate (\ref{apriori}) is used in an essential way in the
computation of the homology of the resulting chain complex. Indeed,
estimate (\ref{apriori}) allows us to prove that this homology is
stable with respect to $C^0$ perturbations of the Lorentzian
metric. Then the conclusion follows from a connection argument,
together with the fact that in the
particular case of a Lorentzian metric of the form
$h=\alpha_{ij}(x)dx^i\otimes dx^j - dy^2$ the resulting
homology is just the homology of the space of based loops on $X$, by
the Morse complex reinterpretation of classical infinite dimensional
Morse theory.  

\section{The Morse complex: an abstract setting}

In this section we give a brief account of the construction of the
Morse complex for a class of gradient-like flows on an infinite
dimensional manifold. The general theory - which holds for a much larger
class - is fully described in \cite{ama05} and \cite{ama08p}. 

Let $\mathscr{X}$ be a Hilbert manifold, endowed with a complete Riemannian
structure, let $\mathscr{Y}$ be an
affine Hilbert space modeled on the Hilbert space $H$, and let 
$\mathscr{M}=\mathscr{X} \times
\mathscr{Y}$ be the product manifold, endowed with the product
Riemannian structure. We denote by $P: T \mathscr{M} \rightarrow H$ 
the map
\[
(x,y,\xi,\eta) \mapsto \eta, \quad x\in \mathscr{X}, \; y\in
\mathscr{Y}, \; \xi\in T_x \mathscr{X}, \; \eta \in T_y \mathscr{Y} =
H,
\]
that is the differential of the projection of $\mathscr{X} \times
\mathscr{Y}$ onto the second factor.

The set of zeroes of a vector field $F$ on $\mathscr{M}$ is 
denoted by $\rest (F)$.
A smooth vector field $F$ on $\mathscr{M}$ is said to be a {\em Morse vector
field} if for every $z\in \rest (F)$ the Jacobian of $F$ at $z$, 
denoted by $\nabla F(z)$, is an
infinitesimally hyperbolic operator on $T_z \mathscr{M}$ (i.e.\ its
spectrum is disjoint from the imaginary axis $i\R$). The positive
and the negative eigenspaces of $\nabla F(z)$ are denoted by
$V^+(\nabla F(z))$ and $V^-(\nabla F(z))$, respectively.  

The set of critical points of a smooth function $f:\mathscr{M}
\rightarrow \R$ is denoted by $\crit(f)$. A smooth function
$f:\mathscr{M} \rightarrow \R$ is said to be a {\em Lyapunov function} for
$F$ if $Df(p)[F(p)]<0$ for every $p$ in $\mathscr{M}$ which is not a
rest point of $F$. If moreover for every critical point $z$ the Hessian
$\nabla^2 f(z)$ is an isomorphism (i.e.\ $f$ is a Morse function), and
$-\nabla^2 f (z)$ is strictly positive (resp.\ strictly negative) on 
$V^+(\nabla F(z))$ (resp.\ on $V^-(\nabla F(z))$), $f$ is said to be a
{\em non-degenerate  Lyapunov function} for $F$. In this case, $\rest
(F) = \crit (f)$, and $F$ is said a (negative) {\em pseudo-gradient}
for $f$. 

A {\em Palais-Smale sequence} for the pair $(F,f)$ is a sequence
$(p_n)\subset \mathscr{M}$ such that $f(p_n)$ is bounded and $Df(p_n)
[F(p_n)]$ is infinitesimal.

Let $\mathscr{A}$ be an open subset of $\mathscr{M}$, and let
$a:\mathscr{M} \rightarrow \R$ be a smooth function such that
$\mathscr{A} = \{a<0\}$. Let $F$ be a smooth Morse vector field on 
$\mathscr{M}$, and let $f$ be a smooth non-degenerate Lyapunov 
function for $F$. We assume the following conditions:

\begin{enumerate}

\item[(F1)] $F$ is bounded on $\mathscr{M}$, and Lipschitz continuous on
  bounded subsets of $\mathscr{M}$. 

\item[(F2)] $F$ has no rest points on the closure of $\mathscr{A}$,
  and $a$ is a Lyapunov function for $F$ on $\mathscr{A}$.

\item[(F3)] $f$ and $a$ are bounded on bounded subsets of $\mathscr{M}$.

\item[(F4)] For every $c,c'\in \R$ the subset $\{f\leq c\} \cap 
\{ a\geq c'\}$ is bounded.

\item[(F5)] Every bounded Palais-Smale sequence for $(F,f)$ has a
  converging subsequence.

\item[(F6)] For every $z=(x,y)\in \rest (F)$, the subspaces 
  $V^-(\nabla F(z))$ and $V^+(\nabla F(z))$ are compact 
  perturbations\footnote{A
  closed linear subspace $V$ of a Hilbert space is said to be a
  compact perturbation of a closed linear subspace $W$ if the
  difference of their orthogonal projectors is a compact operator.} 
  of $T_x \mathscr{M} \times (0)$ and $(0) \times T_y \mathscr{Y} = (0) \times
  H$, respectively.

\item[(F7)] For every $y\in \mathscr{Y}$, the map $\mathscr{X} \rightarrow
  H$, $x\mapsto P F(x,y)$, is compact (i.e.\ it maps bounded sets into
  pre-compact sets).

\end{enumerate}

We denote by $\phi^F: \R \times \mathscr{M} \rightarrow \mathscr{M}$ the
flow of $F$ (which is globally defined because $F$ is bounded and
$\mathscr{M}$ is complete).
Condition (F2) implies that $\mathscr{A}$ is positively invariant for
$\phi^F$. Conditions (F3) and (F4) imply that $f$ is bounded below on 
$\mathscr{M} \setminus \mathscr{A}$, so by (F2) and (F5) $F$ has 
finitely many rest points in $\{f\leq c\}$, for every $c\in \R$. 

By (F6), for every rest point $z=(x,y)$ 
the pair of subspaces $(V^+(\nabla F(z)),T_x \mathscr{X} \times (0))$
is a Fredholm pair\footnote{A pair of closed subspaces $(V,W)$ of a
  Hilbert space $H$ is said to be a Fredholm pair if $V\cap W$ has finite
  dimension, and $V+W$ is closed and has finite codimension. The
  Fredholm index of $(V,W)$ is the integer $\ind (V,W) = \dim V\cap W
  - \dim H/(V+W)$.} in $T_z \mathscr{M}$, and we define  the
{\em relative Morse index} of $z$ as the Fredholm index of such a pair,
\begin{equation}
\label{relmorseind}
i_{\mathrm{rel}} (z):= \ind \bigl(V^+(\nabla F(z)),
T_x \mathscr{X} \times (0)\bigr).
\end{equation}
The stable and unstable manifolds of a rest point $z$, that is the sets
\[
W^s(z;F) := \set{p\in \mathscr{M}}{\lim_{t\rightarrow +\infty}
  \phi^F(t,p)=z}, \quad W^u(z;F) := \set{p\in
  \mathscr{M}}{\lim_{t\rightarrow -\infty} \phi^F(t,p)=z},
\]
are smoothly embedded submanifolds of $\mathscr{M}$, diffeomorphic to
the Hilbert spaces $V^-(\nabla F(z))$ and $V^+(\nabla F(z))$,
respectively (because the Morse vector field $F$ has a Lyapunov
function). Up to a perturbation, we may also assume that the vector
field $F$ has the Morse-Smale property, meaning that the stable and 
unstable manifolds of any two rest points intersect transversally. 
In this situation, $W^u(z;F)\cap W^s(z';F)$ - if non-empty - is a manifold
of dimension $i_{\mathrm{rel}}(z)-i_{\mathrm{rel}}(z')$. 
An orientation of the determinant bundle\footnote{The determinant
  bundle over the space of all Fredholm pairs of a Hilbert space is a
  real line bundle whose fiber at $(V,W)$ is $\Lambda^{\mathrm{top}}
  (V\cap W) \otimes \bigl(\Lambda^{\mathrm{top}} (H/(V+W))\bigr)^*$,
  where $\Lambda^{\mathrm{top}}(X)$ is the degree-$n$ component of the exterior
  algebra of the $n$-dimensional vector space $X$.}
 over the Fredholm pair $(T_z W^s(z;F),(0)\times H)$, for every $z\in \rest
(F)$, induces an orientation of each intersection $W^u(z;F)\cap
W^s(z';F)$.

Every intersection $W^u(z;F)\cap W^s(z';F)$ is pre-compact. 
Indeed, by (F2) this intersection is contained in $\{a\geq
0\}$, and since it is also contained in $\{f\leq f(z)\}$, it is
bounded by (F4). Then
(F5), (F6), and (F7) imply that such an intersection has compact
closure in $\mathscr{M}$. Compactness and transversality imply that,
when $i_{\mathrm{rel}}(z)=i_{\mathrm{rel}}(z')+1$, $W^u(z;F)\cap
W^s(z';F)$ consists of finitely many flow orbits. In this case, we can 
define an integer $d_F(z,z')$ by counting each orbit as $+1$ if its 
orientation agrees with the flow direction, as $-1$ otherwise. 
Let $C_k(F)$ be the free Abelian group generated by the rest points of
$F$ of relative Morse index $k$, and let
$\partial_k : C_k(F) \rightarrow C_{k-1} (F)$ 
be the homomorphism defined by
\[
\partial_k z = \sum_{\substack{z' \in \rest (F) \\
i_{\mathrm{rel}}(z')=k-1}} d_F(z,z')\, z',
\]
for every $z\in \rest (F)$ with $i_{\mathrm{rel}}(z)=k$. The above sum
is finite because $F$ has finitely many rest points in $\{f\leq f(z)\}$. 
These homomorphisms satisfy $\partial_k \circ \partial_{k+1} =0$, so
\[
\{C_k(F),\partial_k\}_{k\in \Z}
\]
is a chain complex of Abelian groups, called the {\em Morse complex of
$F$}. Changing some of the orientations of the
determinant lines over $(T_z W^s(z;F),(0)\times H)$, for $z\in \rest
(F)$, produces an isomorphic chain complex.
Changing the vector field $F$ while keeping the same Lyapunov
function $f$ - hence the same rest points and the same groups $C_k$ - 
produces an isomorphic chain complex (the isomorphisms here being 
non-trivial). Therefore, the homology of the Morse complex depends
only on the Lyapunov function $f$, and it is called
the {\em Morse homology of $f$},
\[
H_k(f) := H_k (\{C_*(F),\partial_*\}) = \ker \partial_k/ \im \partial_{k+1}.
\]

\begin{rem}
\label{montreal}
Assume that the data $F$, $f$, $\mathscr{A}$, $a$ satisfy (F1)-(F7),
and that the function $f$ and the vector field $F$ have the form
\[
f(x,y) = f^+(x) - f^-(y),  \quad F(x,y) = F^+(x) \oplus (-F^-(y)).
\]
where $f^-$, $f^+$ are bounded below, and $f^-$ has finitely many 
critical points. By (F6), all the critical points of $f^+$ and $f^-$
have finite Morse index, and $i_{\mathrm{rel}}(x,y) = i(x) - i(y)$,
$i$ denoting the usual Morse index. In this case it is easily seen that 
\[
C_k(F) = \bigoplus_{p - q = k} C_p(F^+) \otimes C_q(F^-), 
\]
and that the Morse homology of $f$ is 
\[
H_k(f) = \bigoplus_{p - q = k} H_p(f^+) \otimes H_q(f^-).
\]
The Morse homologies of $f^+$ and $f^-$ are isomorphic to the singular 
homologies of $\mathscr{X}$ and $\mathscr{Y}$, respectively (this
holds in general for a Morse function bounded from below, satisfying the
Palais-Smale condition, and having critical points with finite Morse 
index, see \cite{ama06m}). Since $\mathscr{Y}$ has the homology of a
point, we conclude that $H_k(f)$ is isomorphic to $H_k(\mathscr{X})$.
\end{rem}  

We conclude this section by describing the functorial properties of
the Morse complex.
Assume that $(F_0,f_0,\mathscr{A}_0,a_0)$ and
$(F_1,f_1,\mathscr{A}_1,a_1)$ satisfy conditions (F1)-(F7) on the same
Hilbert manifold $\mathscr{M}$. Furthermore, assume the following condition:

\begin{enumerate}

\item[(F8)] For every $c,c'\in \R$ the subset $\{f_0\leq c\} \cap
  \{a_1\geq c'\}$ is bounded.

\end{enumerate}

This condition implies that for
every $z_0\in \rest(F_0)$ and $z_1\in \rest(F_1)$ the intersection
\[
W^u(z_0;F_0) \cap W^s(z_1;F_1)
\] 
is bounded. Indeed, by (F2) this
intersection is contained in $\{a_1\geq 0\}$, and being contained also
in $\{f_0\leq f_0(z_0)\}$, it is bounded by (F8). By (F3) and (F4),
also its positive evolution by
$\phi^{F_1}$ and its negative evolution by $\phi^{F_0}$ are bounded.
Then (F5), (F6), and (F7) imply that
this intersection is pre-compact. Up to perturbing the vector
fields $F_0$ and $F_1$ without affecting their Morse complexes, we may
assume that the above intersections are also transverse. Then (F6) and
(F7) imply that $W^u(z_0;F_0) \cap W^s(z_1;F_1)$ is an oriented
submanifold of dimension
\[
\dim W^u(z_0;F_0) \cap W^s(z_1;F_1) = i_{\mathrm{rel}}(z_0;F_0) - 
i_{\mathrm{rel}}(z_1;F_1).
\]
In particular, when $i_{\mathrm{rel}}(z_0;F_0) =
i_{\mathrm{rel}}(z_1;F_1)$ this intersection is a finite set of
points, each of which carries an orientation sign $+1$ or $-1$. The
sum of these numbers defines an integer $d_{F_0,F_1}(z_0,z_1)$, and
the homomorphism defined by
\[
\Phi_{F_0, F_1} : C_k(F_0) \rightarrow C_k(F_1) , \quad \Phi_{F_0,
  F_1} z_0 = \sum_{\substack{z_1 \in \rest(F_1) \\
    i_{\mathrm{rel}}(z_1;F_1) = k}}d_{F_0,F_1}(z_0,z_1)\, z_1,
\]
for every $z_0\in \rest(F_0)$, $i_{\mathrm{rel}}(z_0;F_0)=k$, 
is a chain map from the Morse complex of $F_0$ to the Morse complex of
$F_1$. Changing the vector fields $F_0$ and $F_1$ while keeping
the same Lyapunov functions $f_0$ and $f_1$, produces homotopic chain
maps, therefore the induced homomorphism between the Morse homology
groups 
\[
\Phi_{f_0,f_1} : H_k(f_0) \rightarrow H_k(f_1)
\]
is well-defined. If $f_0=f_1$, $\Phi_{f_0,f_1}$ is the identity.

Finally, assume that $(F_0,f_0,\mathscr{A}_0,a_0)$,
$(F_1,f_1,\mathscr{A}_1,a_1)$, and  $(F_2,f_2,\mathscr{A}_2,a_2)$,
satisfy (F1)-(F7), and that the pairs $(f_0,a_1)$ and
$(f_1,a_2)$ satisfy (F8). Assume also:

\begin{enumerate}

\item[(F9)] There holds $a_1\leq f_1$, and
for every $c,c'\in \R$ the subsets $\{f_0\leq c\} \cap \{f_1+a_1\geq
c'\}$ and $\{f_1+a_1\leq c\} \cap \{a_2\geq c'\}$ are bounded.

\end{enumerate}

Since
\[
\{f_0\leq c\} \cap \{a_2\geq c'\} \subset \Bigl(\{f_0 \leq c\} \cap
\{f_1+a_1\geq 0\}\Bigr) \cup \Bigl( \{f_1+a_1\leq 0\} \cap \{a_2\geq
c'\} \Bigr),
\]
(F9) implies that also the pair $(f_0,a_2)$ satisfies (F8).
Given $z_0\in \rest(F_0)$ and $z_2\in \rest(F_2)$, we consider the set
\[
W(z_0,z_2) = \set{ (z,z') \in W^u(z_0;F_0) \times W^s(z_2;F_2) }{
  \exists t\geq 0 \mbox{ such that } \phi^{F_1}(t,z)=z' }.
\]
Denote by $W_1$ and $W_2$ the projections of $W(z_0,z_2)$ onto the
first and the second factor of $\mathscr{M} \times \mathscr{M}$.
We claim that $f_1+a_1$ is bounded below on $W_1$. If not, there 
exists a sequence $(w_n,w_n')\subset
W(z_0,z_2)$ such that $f_1(w_n) + a_1(w_n)\rightarrow -\infty$. Since
$a_1\leq f_1$, also $a_1(w_n)\rightarrow -\infty$, so $w_1$ eventually
belongs to $\mathscr{A}_1$. On this set both $f_1$ and $a_1$ are
Lyapunov functions for $F_1$, so also their sum is a Lyapunov function
for $F_1$, and we deduce that $f_1(w_n')+a_1(w_n')\rightarrow
-\infty$. In particular, $f_1(w_n') + a_1(w_n')$ is bounded above. By
condition (F2) for $F_2$, $a_2(w_n')\geq 0$, so (F9) implies that
$(w_n')$ is bounded. But then $f_1(w_n')+a_1(w_n')$ is bounded
(condition (F3)), contradicting the fact that this sequence tends to
$-\infty$. This proves the claim. 

This fact, together with the fact that $f_0$ is bounded above on $W_1$,
implies by (F9) that $W_1$ is bounded. Then $f_1$ is bounded on $W_1$
(condition (F3)), hence it is bounded above on $W_2$. 
Since $a_2$ is bounded below on $W_2$, condition (F8) for the pair
$(f_1,a_2)$ implies that also $W_2$ is bounded. By (F3) and (F4), also
the negative evolution of $W_1$ by $\phi^{F_0}$, the positive evolution of
$W_2$ by $\phi^{F_2}$, and the set of orbits of $\phi^{F_1}$
connecting $W_1$ to $W_2$ are bounded.

As before, conditions (F5), (F6), and (F7), imply that $W(z_0,z_2)$ has
compact closure in $\mathscr{M}\times \mathscr{M}$, and up to 
perturbations, it is an oriented submanifold of dimension
\[
\dim W(z_0,z_2) = i_{\mathrm{rel}}(z_0;F_0) -
i_{\mathrm{rel}}(z_2;F_2) +1.
\]
When $i_{\mathrm{rel}}(z_0;F_0) = i_{\mathrm{rel}}(z_2;F_2) - 1$,
$W(z_0,z_2)$ is a finite set of points with orientation signs, and
their sum defines an integer $d_{F_0,F_1,F_2}(z_0,z_2)$. The
homomorphism
\[
\Psi : C_k(F_0) \rightarrow C_{k+1}(F_2), \quad \Psi z_0 = \sum_{\substack{z_2
  \in \rest(F_2) \\ i_{\mathrm{rel}}(z_2;F_2)=k+1}}
d_{F_0,F_1,F_2}(z_0,z_2) \, z_2,
\]
is a chain homotopy between $\Phi_{F_0,F_2}$ and $\Phi_{F_1,F_2} \circ
\Phi_{F_0, F_1}$. Therefore, 
\[
\Phi_{f_1,f_2} \circ \Phi_{f_0,f_1} = \Phi_{f_0,f_2} : H_k(f_0)
\rightarrow H_k(f_2),
\]
for every $k\in \Z$.

\begin{rem}
Condition (F9) can be replaced by the more symmetric assumption that
the set $\{a_1\leq
c\} \cap \{a_2\geq c'\}$ should be bounded for every $c,c'\in \R$. For the
purposes of this paper however (F9) is easier to check. 
\end{rem}

\section{The main result}

Let $(X,g)$ be a smooth connected complete Riemannian manifold, and let 
$M=X\times \R$ be endowed with the Lorentzian structure
\[
h(x,y)[(\xi_1,\eta_1),(\xi_2,\eta_2)] =
g(x)[\alpha(x,y)\xi_1,\xi_2] - \beta(x,y) \eta_1 \eta_2,
\]
where $(x,y)\in X \times \R$, $\xi_1,\xi_2\in T_x X$, 
$\eta_1,\eta_2\in T_y\R = \R$, 
$\alpha$ is a smooth section of the bundle over
$X\times \R$ of $g$-symmetric and $g$-positive $(1,1)$ tensors on $X$,
and $\beta$ is a smooth positive real function on $M$. 
In order to emphasize
the role of $\alpha$ and $\beta$, the Lorentzian 
energy functional is denoted by
\[
E_{\alpha,\beta}(x,y) = \frac{1}{2} \int_0^1 g(x)[\alpha(x,y)x',x']\,
dt - \frac{1}{2} \int_0^1 \beta(x,y)|y'|^2\, dt,
\]
for a curve $(x,y):[0,1]\rightarrow M$.
The aim of this paper is to prove the following:

\begin{thm}
\label{main}
Let $(X,g)$ be a smooth connected complete Riemannian manifold, and let 
$M=X\times \R$ be endowed with the Lorentzian structure $h$ defined as
above. Assume that:
\begin{enumerate}
\item[(a0)] there is a positive number $s_0$ such that $X \times ]-s,s[$ is
  convex for every $s>s_0$;
\item[(a1)] there are positive numbers
  $\underline{\alpha},\overline{\alpha}$ such that
  $\underline{\alpha} I \leq \alpha \leq \overline{\alpha} I$ as
  $g(x)$-symmetric $(1,1)$ tensors, on $X \times [-s_0-1,s_0+1]$;
\item[(a2)] there are positive numbers
  $\underline{\beta},\overline{\beta}$ such that $\underline{\beta}
  \leq \beta \leq \overline{\beta}$ on $X \times [-s_0-1,s_0+1]$;
\item[(a3)] the function $\xi \mapsto g(\partial_y \alpha\, \xi,\xi)/g(\alpha\,
  \xi,\xi)$ is bounded on $(TX\setminus \{0\})\times [-s_0-1,s_0+1]$;
\item[(a4)] the function $\partial_y \beta/\beta$ is bounded on $X \times
  [-s_0-1,s_0+1]$.
\end{enumerate}
Let $z_0=(x_0,y_0)$ and $z_1=(x_1,y_1)$ be two points in $X \times
]-s_0,s_0[$, assumed to be non-conjugate along every geodesic joining
them. Given $k\in \Z$, let $C_k$ be the free Abelian group generated 
by the geodesics $\gamma$ joining $z_0$ and $z_1$ of relative Morse 
index $k$. Then there is a pseudo-gradient vector field $F$ for
$E_{\alpha,\beta}$ on the Hilbert manifold of $W^{1,2}$ curves joining
$z_0$ to $z_1$ such that $C_k(F)=C_k$, and the homology of the
corresponding Morse complex $\{C_k(F),\partial_k\}_{k\in \Z}$ is isomorphic
to the singular homology of $\Omega(X)$, the based loop space of $X$.
\end{thm}

Notice that assumptions (a1)-(a4) always hold if $X$ is
compact, so Theorem \ref{main} implies the theorem stated in
the introduction. 
By condition (a0), all the geodesics joining $z_0$ and $z_1$ are
contained in $X \times [-s_0,s_0]$.

If $\gamma=(x,y):[0,1] \rightarrow M$ is a geodesic, the
Euler-Lagrange equation for $E_{\alpha,\beta}$ with respect to 
variations on the second component of $(x,y)$ is the equation
\begin{equation} 
\label{equa}
(\beta(x,y)y')' + \frac{1}{2} g(x)[\partial_y \alpha (x,y) x',x'] - \frac{1}{2}
\partial_y \beta(x,y)|y'|^2 =0.
\end{equation}
Moreover, the conservation of energy produces the identity
\begin{equation}
\label{energy}
\frac{1}{2} g(x(t))[\alpha(x(t),y(t))x'(t),x'(t)] - \frac{1}{2}
\beta(x(t),y(t)) |y'(t)|^2 = E_{\alpha,\beta}(x,y), \quad \forall
t\in [0,1].
\end{equation}

The convexity assumption (a0) can be restated in terms of the
derivatives of the Lorentzian structure:

\begin{lem}
\label{conv}
The sets $X \times ]-s,s[$ are convex for every $s>s_0$ if
and only if for every $x\in X$
\begin{equation}
\label{segno}
\partial_y \alpha(x,y) \leq 0 \quad \forall y>s_0, \quad \partial_y
\alpha(x,y) \geq 0 \quad \forall y<-s_0,
\end{equation}
as $g(x)$-symmetric $(1,1)$ tensors.
\end{lem}

\begin{proof}
First notice that $X \times ]-s,s[$ is convex if and only if every geodesic
$\gamma=(x,y):]-\epsilon,\epsilon[ \rightarrow X \times \R$ with
$|y(0)|=s$ and $|y(t)|\leq s$ for every $t$ satisfies $|y(t)|=s$ for
every $t$. 

Assume that $\partial_y \alpha(x,y)\leq 0$ does not hold at some point
$(x_0,y_0)\in X \times \R$, that is there is some $\xi\in T_{x_0}X$
such that
\[
g(x_0) [\partial_y \alpha(x_0,y_0)\xi,\xi]>0.
\]
Let $(x,y):]-\epsilon,\epsilon[ \rightarrow X \times \R$ be the
geodesic such that $(x(0),y(0))=(x_0,y_0)$ and
$(x'(0),y'(0))=(\xi,0)$. The geodesic equation (\ref{equa})
for $y$ at $t=0$ yields 
\[
\beta(x_0,y_0) y''(0) +
\frac{1}{2} g(x_0)[\partial_y \alpha (x_0,y_0)\xi,\xi] = 0.
\]
Therefore $y''(0)<0$, hence $y_0$ is a strict local maximum of $y$,
and the set $X \times ]-\infty,y_0[$ is not convex. So, if $X \times ]-s,s[$
is convex then $\partial_y \alpha(x,s)\leq 0$ for every $x\in
X$. Similarly, $\partial_y \alpha(x,-s)\geq 0$ for every $x\in X$.  

To prove the converse, assume (\ref{segno}) and choose some
$s>s_0$. Let $(x,y)$ be a geodesic in $X\times \R$ with $y(0)=s$ and
$y(t)\leq s$ for every $t$ (the case $y(0)=-s$ and $y(t)\geq -s$ for
every $t$ being analogous). There is an interval $I$, neighborhood of
$0$, on which $g(x)[\partial_y \alpha(x,y) x',x']\leq 0$. From the
geodesic equation (\ref{equa}) we deduce
\[
y''(t)+ u(t) y'(t) \geq 0, \quad \forall y\in I,
\]
for a suitable smooth function $u$. Since $y(t)\leq y(0)=s$ for every
$t\in I$, the strong maximum principle implies that $y(t)=s$ for every
$t\in I$. This shows that $X \times ]-s,s[$ is convex.
\end{proof} \qed

Up to changing the Lorentzian structure outside $X\times [-s_0,s_0]$
without affecting the convexity assumption (a0), we may assume that the
bounds listed in (a1)-(a4) hold on the whole $X\times \R$. Indeed, it
is enough to choose a smooth function $\psi: \R \rightarrow \R$ such
that $\psi=1$ on $[-s_0,s_0]$, $\psi=0$ on $\R \setminus
[-s_0-1,s_0+1]$, $\psi'\geq 0$ on $[-s_0-1,-s_0]$, $\psi' \leq 0$ on
$[s_0,s_0+1]$, and to replace $\alpha$ and $\beta$ by
\[
\psi(y) \alpha(x,y) + (1-\psi(y)) \underline{\alpha} I,\quad \mbox{and} \quad  
\psi(y) \beta(x,y) + (1-\psi(y)) \underline{\beta}.
\]
Conditions (a1)-(a4) now hold on $X\times \R$. 
The derivative with respect to $y$ of the second function is
\[
\psi \partial_y \alpha + \psi' (\alpha-\overline{\alpha}), 
\]
so the new Lorentzian structure satisfies also (a0), by Lemma
\ref{conv}.  

By the above considerations, we may assume that the coefficients of
our Lorentzian structure belong to the set 
$\Gamma = \Gamma(s_0,\underline\alpha,\overline{\alpha},\underline{\beta},
\overline{\beta},a,b)$ consisting of all $(\alpha,\beta)$ satisfying:

\begin{enumerate}

\item[(h0)] $\partial_y \alpha(x,y)\leq 0$ for every $y\geq s_0$,
  $\partial_y \alpha (x,y)\geq 0$ for every $y\leq -s_0$, for every
  $x\in X$;

\item[(h1)] $\underline{\alpha} I \leq \alpha \leq \overline{\alpha}
  I$ on $M$;

\item[(h2)] $\underline{\beta} \leq \beta \leq \overline{\beta}$ on $M$;

\item[(h3)] $\displaystyle{\sup_{\xi \in T 
      X \setminus \{0\}} } \frac{|g(\partial_y \alpha \, \xi,\xi)|}{
  g(\alpha \, \xi,\xi)} \leq a$ on $M$;

\item[(h4)] $|\partial_y \beta|/\beta \leq b$ on $M$.

\end{enumerate}

\section{A priori bounds}

The following a priori estimate sharpens Theorem 4.1 in \cite{bfm94}:

\begin{lem}
\label{stima}
Let $(\alpha,\beta)\in 
\Gamma(s_0,\underline{\alpha},\overline{\alpha},\underline{\beta},
\overline{\beta},a,b)$ and let $y_0,y_1\in
]-s_0,s_0[$. Then there are constants
\begin{eqnarray*}
p_1 = p_1(s_0,\underline{\beta},\overline{\beta},a,b) 
= \frac{a}{a+b} \frac{\overline{\beta}}{\underline{\beta}^2}
e^{2(a+b)s_0}, \\ q_1 =
q_1(s_0,\underline{\beta},\overline{\beta},a,b,
y_1-y_0) = 
\frac{\overline{\beta}^2}{\underline{\beta}^2} |y_1-y_0|^2
  e^{2(a+b)s_0},
\end{eqnarray*}
such that every geodesic $\gamma=(x,y):[0,1]\rightarrow M$ with
$y(0)=y_0$ and $y(1)=y_1$ and $E_{\alpha,\beta}(x,y)=c$ satisfies
\[
\|y\|_{\infty} \leq s_0, \quad \|y'\|_{\infty}^2 \leq p_1
c^+ + q_1, \quad c \geq - \frac{1}{2}
\overline{\beta} q_1,
\]
where $c^+$ denotes the positive part of the real number $c$.
\end{lem}

\begin{proof}
Since $X\times ]-s,s[$ is convex for every $s>s_0$ and $|y_0|<s_0$,
$|y_1|<s_0$, we have
\begin{equation}
\label{uni}
\|y\|_{\infty} \leq s_0.
\end{equation}
Let $t_1\in [0,1]$ be a point where $|y'|$ attains its maximum. 
By the mean value theorem, up to
inverting the time parameterization, $t\mapsto 1-t$, we may assume
that there is $t_0\leq t_1$ such that $y'$ does not change sign on 
$[t_0,t_1]$, and $y'(t_0)$ is either $0$ or $y_1-y_0$. Up to
considering the change of variable $(x,y)\mapsto (x,-y)$, we may
also assume that $y'\geq 0$ on $[t_0,t_1]$.

By the geodesic equation (\ref{equa}), by (h3), (h4), and by
the conservation of energy (\ref{energy}), we have
\begin{eqnarray*} 
(\beta y')' = \frac{1}{2} \partial_y \beta |y'|^2 - \frac{1}{2}
g[\partial_y \alpha x',x'] \leq \frac{b}{2} \beta |y'|^2 +
\frac{a}{2} g[\alpha x',x'] \\ = \frac{b}{2}  \beta |y'|^2 +
\frac{a}{2}  \beta |y'|^2 + \frac{a}{2}  c \leq 
\frac{a+b}{2} \beta |y'|^2 + \frac{a}{2} c^+,
\end{eqnarray*}
where $c^+$ denotes the positive part of the real number $c$. Then
\[
(\beta y')' - \frac{a+b}{2} \beta |y'|^2 \leq \frac{a}{2} c^+,
\]
and multiplying both sides of this inequality by $2\beta y'
e^{-(a+b)y}$ we obtain, using also (h2),
\[
\left( \beta^2 |y'|^2 e^{-(a+b)y} \right)' \leq \frac{a}{a+b} c^+ \beta
\left(-e^{-(a+b)y} \right)' \leq \frac{a}{a+b} c^+ \overline{\beta}
\left(-e^{-(a+b)y} \right)',
\]
on the interval $[t_0,t_1]$. Since $y'(t_0)\leq |y_1-y_0|$,
integration on such an interval gives us
\begin{eqnarray*}
\beta(x(t_1),y(t_1))^2 |y'(t_1)|^2 e^{-(a+b)y(t_1)} \\ \leq
\beta(x(t_0),y(t_0))^2 |y_1-y_0|^2 e^{-(a+b)y(t_0)} + \frac{a}{a+b}
c^+ \overline{\beta} \left( e^{-(a+b)y(t_0)} - e^{-(a+b)y(t_1)}
\right).
\end{eqnarray*}
By (h2) and (\ref{uni}), we deduce that
\[
\|y'\|_{\infty}^2 = |y'(t_1)|^2 \leq \left( \frac{a}{a+b}
  \frac{\overline{\beta}}{\underline{\beta}^2} c^+ +
  \frac{\overline{\beta}^2}{\underline{\beta}^2} |y_1-y_0|^2 \right)
e^{2(a+b)s_0},
\]
hence
\[
\|y'\|_{\infty}^2 \leq p_1 c^+ + q_1,
\]
as claimed. In particular, if $c\leq 0$ we have $\|y'\|_{\infty}^2 
\leq q_1$. By (h2) this implies
\[
c = E_{\alpha,\beta}(x,y) \geq - \frac{1}{2} \overline{\beta}
\|y'\|_2^2 \geq  - \frac{1}{2} \overline{\beta}
\|y'\|_{\infty}^2 \geq  - \frac{1}{2} \overline{\beta} q_1,
\]
concluding the proof.
\end{proof} \qed

\begin{rem}
\label{stmar}
By Lemma \ref{conv}, if the Lorentzian structure defined by
$(\alpha,\beta)$ satisfies (h0), so does the Lorentzian structure defined by
$(\alpha,\beta+\lambda)$, for every positive constant
$\lambda$. Therefore, if $(\alpha,\beta) \in
\Gamma(s_0,\underline{\alpha},\overline{\alpha},\underline{\beta}, 
\overline{\beta},a,b)$ then $(\alpha,\beta+\lambda) \in  
\Gamma(s_0,\underline{\alpha},\overline{\alpha},
\underline{\beta}+\lambda,\overline{\beta}+\lambda,a,b)$, for every
$\lambda \geq 0$.  
Notice also that since
$\overline{\beta} \geq \underline{\beta}$, the inequalities  
\[
\frac{\overline{\beta} + \lambda}{\underline{\beta} + \lambda} 
= 1 +
\frac{\overline{\beta}-\underline{\beta}}{\underline{\beta}+\lambda}
  \leq 1 + \frac{\overline{\beta}-
    \underline{\beta}}{\underline{\beta}} =
\frac{\overline{\beta}}{\underline{\beta}}, \quad
\frac{\overline{\beta} + \lambda}{(\underline{\beta} + \lambda)^2}  \leq
\frac{\overline{\beta}}{\underline{\beta}} \frac{1}{\underline{\beta}
  + \lambda} \leq 
\frac{\overline{\beta}}{\underline{\beta}^2},  
\]
hold for any $\lambda\geq 0$. We conclude that
\begin{eqnarray*}
p_1(s_0,\underline{\beta}+\lambda,\overline{\beta}+\lambda,a,b) \leq 
p_1(s_0,\underline{\beta},\overline{\beta},a,b), \\
q_1(s_0,\underline{\beta}+\lambda,\overline{\beta}+\lambda,a,b,y_1-y_0) 
\leq
q_1(s_0,\underline{\beta},\overline{\beta},a,b,y_1-y_0),
\end{eqnarray*}
for every positive constant $\lambda$.
\end{rem} 
 
Lemma \ref{stima} has the following consequence:

\begin{lem}
\label{stima2}
Let $(\alpha,\beta) \in \Gamma(s_0,\underline{\alpha},
\overline{\alpha},\underline{\beta},\overline{\beta},a,b)$, and 
let $y_0,y_1\in ]-s_0,s_0[$.
For every $\lambda$ in the interval
$[0,2/p_1(s_0,\underline{\beta},\overline{\beta},a,b)]$, 
for every $\mu\in [0,\lambda]$, and for every critical point $(x,y)$ of 
$E_{\alpha,\beta+\mu}$, there holds
\[
E_{\alpha,\beta+\lambda}(x,y) \geq
  -\frac{1}{2} (\overline{\beta}+\lambda)
  q_1(s_0,\underline{\beta},\overline{\beta},a,b,y_1-y_0).
\]
\end{lem}

\begin{proof}
Let $(x,y)$ be a critical point of $E_{\alpha,\beta+\mu}$, and notice
that
\begin{equation}
\label{s1}
E_{\alpha,\beta+\lambda}(x,y) = E_{\alpha,\beta+\mu}(x,y) -
\frac{1}{2} (\lambda-\mu) \|y'\|_2^2.
\end{equation}
If $E_{\alpha,\beta+\mu}(x,y)\leq 0$, by Lemma \ref{stima} we have
\begin{equation}
\label{s2}
E_{\alpha,\beta+\mu}(x,y) - \frac{1}{2} (\lambda-\mu) \|y'\|_2^2 \geq
- \frac{1}{2}(\overline{\beta}+\mu)q_1(\mu) - 
\frac{1}{2} (\lambda-\mu) q_1(\mu) = -\frac{1}{2}
(\overline{\beta} + \lambda) q_1(\mu),
\end{equation}
where $q_1(\mu)$ stays for
$q_1(s_0,\underline{\beta}+\mu,\overline{\beta}+\mu, a,b,y_1-y_0)$.
By Remark \ref{stmar}, $q_1(\mu)\leq q_1(0)$, so (\ref{s1}) and 
(\ref{s2}) imply
\begin{equation}
\label{s3}
E_{\alpha,\beta+\lambda}(x,y) \geq - \frac{1}{2}
(\overline{\beta}+\lambda) q_1(0).
\end{equation}
If $E_{\alpha,\beta+\mu}(x,y)>0$, by Lemma \ref{stima} we have
\begin{equation}
\label{s4}
\begin{split}
E_{\alpha,\beta+\mu}(x,y) - \frac{1}{2} (\lambda-\mu) \|y'\|_2^2
\geq E_{\alpha,\beta+\mu}(x,y) - \frac{1}{2} (\lambda-\mu)
p_1(\mu) E_{\alpha,\beta+\mu}(x,y) \\ - \frac{1}{2} (\lambda-\mu)
q_1(\mu) \geq  \left(1 - \frac{1}{2} (\lambda-\mu)
  p_1(\mu) \right) E_{\alpha,\beta+\mu}(x,y) -
\frac{1}{2} \lambda q_1(0),
\end{split}
\end{equation}
where $p_1(\mu)$ stays for
$p_1(s_0,\underline{\beta}+\mu,\overline{\beta}+\mu, a,b)$ and
we have used again the inequality $q_1(\mu)\leq
q_1(0)$. By Remark \ref{stmar} we also have 
$p_1(\mu)\leq p_1(0)$, so
\[
\frac{1}{2} (\lambda-\mu) p_1(\mu) \leq \frac{1}{2}
\lambda p_1(0) \leq 1,
\]
by the assumption on $\lambda$.
Therefore, (\ref{s1}) and (\ref{s4}) imply
\[
E_{\alpha,\beta+\lambda}(x,y) \geq - \frac{1}{2} \lambda
q_1(0) \geq - \frac{1}{2} (\overline{\beta}+\lambda) q_1(0).
\]
We conclude that every critical point
$(x,y)$ of $E_{\alpha,\beta+\mu}$ satisfies (\ref{s3}), as claimed.
\end{proof} \qed

\section{The $\mathbf{W^{1,2}}$ setting for the energy functional}

Let $\mathscr{X}$ be the Hilbert manifold of paths $x:[0,1]
\rightarrow X$ of Sobolev class $W^{1,2}$, such that 
$x(0)=x_0$ and $x(1)=x_1$. Let $\mathscr{Y}$ be the affine Hilbert manifold
\[
\mathscr{Y} = \set{y\in W^{1,2}([0,1],\R)}{y(0)=y_0, \; y(1)=y_1},
\]
modeled on the Hilbert space $H=W^{1,2}_0([0,1],\R)$.
The space $\mathscr{M}$ of $W^{1,2}$ paths in $M$ joining $(x_0,y_0)$
to $(x_1,y_1)$ is then the product manifold $\mathscr{M} =
\mathscr{X} \times \mathscr{Y}$. The energy functional $E_{\alpha,\beta}$ is
smooth on $\mathscr{M}$.

On $\mathscr{M}=\mathscr{X} \times \mathscr{Y}$ we consider the Riemannian
structure given by
\begin{equation}
\label{metric}
\langle (\xi_1,\eta_1), (\xi_2,\eta_2) \rangle := \int_0^1
g(\nabla_t^g \xi_1,\nabla_t^g \xi_2)\, dt + \int_0^1 \eta_1' \eta_2'\, dt,
\end{equation}
where $\xi_1,\xi_2 \in T_x \mathscr{X}$ are $W^{1,2}_0$ vector fields
along $x\in \mathscr{X}$, $\nabla_t^g$ denotes the Levi-Civita
covariant derivative along $x$ induced by the metric $g$, and 
$\eta_1,\eta_2\in H$. The associated norm is denoted by $\|\cdot\|$,
and the $\langle\cdot,\cdot\rangle$-gradient of a functional
$f:\mathscr{M}\rightarrow \R$ is denoted by $\nabla f$. 
The Riemannian metric (\ref{metric}) induces a complete distance on
$\mathscr{M}$, and
a subset $\mathscr{A}$ of $\mathscr{M}$ is bounded if and only if
there is some number $c$ such that
\[
\int_0^1 g(x',x')\, dt \leq c, \quad \|y'\|_2 \leq c,
\]
for every $(x,y)\in \mathscr{A}$, $\|\cdot\|_2$ denoting the $L^2$ norm.   

\begin{lem}
\label{exc1}
Let $z=(x,y)$ be a critical point of $E_{\alpha,\beta}$. Then the
positive eigenspace (resp.\ the negative eigenspace) of the Hessian of
$E_{\alpha,\beta}$ at 
$z$ is a compact perturbation of $T_x \mathscr{X} \times (0)$
(resp.\ of $(0) \times T_y \mathscr{Y} = (0) \times H$).
\end{lem}

\begin{proof}
The second differential of the energy functional $E_{\alpha,\beta}$ at
a critical point $z=(x,y)$ is given by the expression
\[
D^2 E_{\alpha,\beta}(x,y)[\zeta_1,\zeta_2] = \int_0^1 \bigl[
  h(\nabla_t^h \zeta_1, \nabla_t^h \zeta_2) -
  h(R(\zeta_1,z')z',\zeta_2) \bigr]\, dt,
\]
for every $\zeta_1,\zeta_2\in T_z \mathscr{M}$, where $\nabla_t^h$
denotes the Levi-Civita covariant derivative along the curve $z$
associated to the Lorentzian structure $h$, and $R$ is the Riemann
tensor of $h$ (see for instance \cite{one83}, Chapter 10, Proposition 39). 
Since any two covariant derivatives along a fixed curve
$z$ differ by a tensor, for every $\zeta=(\xi,\eta)\in T_x \mathscr{X}
\times H = T_z \mathscr{M}$ we have
\[
\nabla_t^h \zeta (t) = (\nabla_t^g \xi(t) + \Gamma_X(t)[\zeta(t)],
\eta'(t) + \Gamma_Y(t)[\zeta(t)]), \quad \forall t\in [0,1],
\]
for some $(1,1)$-tensor $\Gamma_X \oplus \Gamma_Y$ on $z^*(TM)$. Therefore,
writing $\zeta_1=(\xi_1,\eta_1)$, $\zeta_2=(\xi_2,\eta_2)$, with
$\xi_1,\xi_2\in T_x \mathscr{X}$, $\eta_1,\eta_2\in H$, we have
\begin{eqnarray*}
D^2 E_{\alpha,\beta}(z)[\zeta_1,\zeta_2] = \int_0^1 \bigl[
  g(\alpha(z)\nabla_t^g \xi_1, \nabla_t^g \xi_2) - \beta(z)
  \eta_1'\eta_2'\bigr] \, dt \\
+ \int_0^1 \bigl[ g(\alpha(z) \nabla_t^g \xi_1,\Gamma_X [\zeta_2])
+g(\alpha(z) \Gamma_X [\zeta_1], \nabla_t^g \xi_2) + g(\alpha(z)
\Gamma_X [\zeta_1], \Gamma_X [\zeta_2]) \\
- \beta(z) \eta_1' \Gamma_Y[\zeta_2] - \beta(z)
\Gamma_Y[\zeta_1]\eta_2' -  \beta(z) \Gamma_Y[\zeta_1]
\Gamma_Y[\zeta_2] - h(R(\zeta_1,z')z',\zeta_2) \bigr]\, dt.
\end{eqnarray*}
The first integral in the above formula defines a symmetric continuous
bilinear form on $T_z \mathscr{M}$, represented by an invertible
self-adjoint operator $T\in L(T_z \mathscr{M}, T_z \mathscr{M})$,
whose positive and negative eigenspaces are $T_x \mathscr{X} \times
(0)$ and $(0) \times H$, respectively. By the compactness of the
embedding $W^{1,2} \hookrightarrow L^2$, the second integral in the
above formula defines a symmetric bilinear form which is continuous
with respect to the weak topology of $T_z \mathscr{M}$, hence it is
represented by a compact self-adjoint operator $K\in L(T_z
\mathscr{M}, T_z \mathscr{M})$. 

Then $\nabla^2 E_{\alpha,\beta}(z) = T+K$, and the conclusion follows
from the fact that the positive (resp.\ negative) eigenspace of a
compact perturbation of a self-adjoint invertible operator $T$ is a
compact perturbation of the positive (resp.\ negative) eigenspace of
$T$ (see \cite{ama01}, Proposition 2.2).  
\end{proof} \qed 

Therefore, the relative Morse index 
\[
i_{\mathrm{rel}}(z) = \ind (V^-(\nabla^2 E_{\alpha,\beta}(z)), T_x
\mathscr{X} \times (0))
\]
introduced in (\ref{relmorseind}) is a well-defined integer. 

\begin{lem}
\label{interliv}
Let $(\alpha_0,\beta_0)$ and $(\alpha_1,\beta_1)$ be Lorentzian
structures in $\Gamma(s_0,\underline{\alpha},\overline{\alpha},
\underline{\beta},\overline{\beta},a,b)$ such that
\[
\max \{ \|\alpha_0-\beta_0\|_{\infty} , \|\beta_1 - \beta_0\|_{\infty}
\} < \frac{\lambda \underline{\alpha}}{\underline{\alpha} +
  \overline{\beta}},
\]
for some $\lambda>0$. Then for every $c_0,c_1\in \R$ the set 
\[
\{E_{\alpha_0,\beta_0} \leq c_0\} \cap \{ E_{\alpha_1,\beta_1+\lambda}
\geq c_1\}
\]
is bounded in $\mathscr{M}$.
\end{lem}

\begin{proof}
Let $(x,y) \in \{E_{\alpha_0,\beta_0} \leq c_0\} \cap
\{E_{\alpha_1,\beta_1+\lambda}\geq c_1\}$. Then
\begin{equation}
\label{xp}
\int_0^1 g(x',x')\, dt \leq \frac{1}{\underline{\alpha}} \int_0^1
g(\alpha_0 x',x')\, dt = \frac{2}{\underline{\alpha}}
E_{\alpha_0,\beta_0} (x,y) + \frac{1}{\underline{\alpha}} \int_0^1
\beta_0 |y'|^2\, dt \leq \frac{2}{\underline{\alpha}} c_0 +
\frac{\overline{\beta}}{\underline{\alpha}} \|y'\|_2^2.
\end{equation}
Therefore,
\begin{eqnarray*}
2c_1 \leq 2 E_{\alpha_1,\beta_1+\lambda}(x,y) = \int_0^1 g(\alpha_1
x',x')\, dt - \int_0^1 \beta_1 |y'|^2\, dt - \lambda \|y'\|_2^2 \\
\leq 2 E_{\alpha_0,\beta_0} (x,y) + \|\alpha_1-\alpha_0\|_{\infty}
\int_0^1 g(x',x')\, dt - (\lambda - \|\beta_1-\beta_0\|_{\infty} )
\|y'\|_2^2 \\ \leq 2c_0 \left( 1 +
\frac{\|\alpha_1-\alpha_0\|_{\infty}}{\underline{\alpha}} \right) -
\left(\lambda - \|\beta_1-\beta_0\|_{\infty} -
\frac{\overline{\beta}}{\underline{\alpha}}
\|\alpha_1-\alpha_0\|_{\infty} \right) \|y'\|_2^2.
\end{eqnarray*}
By our assumption on $\|\alpha_1-\alpha_0\|_{\infty}$ and
$\|\beta_1-\beta_0\|_{\infty}$, the coefficient of $\|y'\|_2^2$ in the
latter expression is positive, so the above inequality implies that
$\|y'\|_2$ is uniformly bounded. By (\ref{xp}), also $\int_0^1
g(x',x')\, dt$ is uniformly bounded, hence $(x,y)$ is uniformly
bounded in $\mathscr{M}$.
\end{proof} \qed 

The following result is well known (see for instance \cite{bf94b},
Lemma 3.2), but we include a proof for sake of completeness:

\begin{lem}
\label{psobs}
Assume that $\alpha$ and $\beta$ are strictly positive (the first
object as a section over $M$ 
of $g$-symmetric (1,1) tensors on $X$, the second as a real function
on $M$). Then every
bounded sequence $(z_h)\subset \mathscr{M}$ such that $\|\nabla
E_{\alpha,\beta}(z_h)\|\rightarrow 0$ has a converging subsequence.
\end{lem}

\begin{proof}
Let $(z_h)=(x_h,y_h)$ be a bounded sequence in $\mathscr{M}$. By the
compactness of the embedding of $W^{1,2}$ into the space of continuous
curves, up to a subsequence we may assume that $(z_h)$ converges
uniformly to some $\overline{z}=(\overline{x},\overline{y}) \in
\mathscr{M}$. Let $U$ be an open neighborhood of $0$ in $\R^n$, with
$n=\dim X$, and let 
\[
[0,1]\times U \rightarrow [0,1] \times X, \quad (t,\xi) \mapsto
(t,\varphi(t,\xi)), 
\]
be a smooth coordinate system such that $\varphi(0,0)=x_0$,
$\varphi(1,0)=x_1$, and $\overline{x}(t)=\varphi(t,\overline{\xi}(t))$
for every $t\in [0,1]$, for some $\overline{\xi}\in
W^{1,2}_0([0,1],U)$. For instance, such a diffeomorphism can be
constructed by choosing a smooth curve $x: [0,1]\rightarrow X$
connecting $x_0$ and $x_1$, with
$\dist(x(t),\overline{x}(t)) < \rho$ for every $t\in [0,1]$, $\rho$
being the injectivity radius of a large compact neighborhood of
$\overline{x}([0,1])$ in $X$, and by setting
\[
\varphi(t,\xi) := \exp_{x(t)} \bigl(\Phi(t) \xi\bigr),
\]
where $\xi\in \R^n$, $|\xi|<\rho$, and $\Phi$ is a smooth orthogonal
trivialization of the vector bundle $x^*(TX)$ over $[0,1]$.

Then the map
\[
\varphi_* : W^{1,2}_0([0,1],U) \rightarrow \mathscr{X}, \quad
\varphi_*(\xi) (t) := \varphi(t,\xi(t)),
\]
is a smooth local coordinate system on the Hilbert manifold
$\mathscr{X}$. Up to replacing $U$ by a smaller neighborhood, we may
assume that $D\varphi_*$ is bounded together with its inverse from the
standard metric of $W^{1,2}_0([0,1],\R^n)$ to the Riemannian metric of
$\mathscr{X}$ defined in (\ref{metric}). 

Since $(x_h)$ converges uniformly to $\overline{x}$, it eventually
belongs to the image of $\varphi_*$, so up to a subsequence we may
assume that $\varphi_*(\xi_h)=x_h$ for
some $(\xi_h)\subset W^{1,2}_0([0,1],U)$ converging to
$\overline{\xi}$ uniformly and weakly in $W^{1,2}_0([0,1],\R^n)$.
By the properties of $\varphi_*$, the fact that $\|\nabla
E_{\alpha,\beta}(z_n)\|$ is infinitesimal is equivalent to the fact
that, setting
\[
\tilde{E}(\xi,y) := E_{\alpha,\beta}(\varphi_*(\xi),y), \quad \forall \xi\in
W^{1,2}_0([0,1],U), \; \forall y\in \mathscr{Y},
\]
the sequence $(D\tilde{E}(\xi_h,y_h))$ is infinitesimal in the dual of
$W^{1,2}_0([0,1],\R^n) \times W^{1,2}_0([0,1],\R)$. It is sufficient
to show that $(\xi_h,y_h)$ converges to
$(\overline{\xi},\overline{y})$ strongly in $W^{1,2}$,
because this fact is equivalent to the convergence of $(x_h,y_h)$ to
$(\overline{x},\overline{y})$ in $\mathscr{M}$.

A direct computation shows that $\tilde{E}$ has the form
\[
\tilde{E}(\xi,y) = \frac{1}{2} \int_0^1 \left[ A(t,\xi,y)\xi' \cdot
  \xi' + B(t,\xi,y) \cdot \xi' + C(t,\xi,y) \right]\, dt - \frac{1}{2}
  \int_0^1 \tilde{\beta}(t,\xi,y) |y'|^2\, dt,
\]
where $A$ takes value in the space of positive definite symmetric
$n\times n$ matrices, and $\tilde{\beta}$ is strictly positive. Since
$(\xi_h-\overline{\xi})$ is bounded in $W^{1,2}_0$, the sequence
$D\tilde{E}(\xi_h,y_h)[(\xi_h-\overline{\xi},0)]$ is
infinitesimal. From the expression
\begin{eqnarray*}
D\tilde{E}(\xi_h,y_h)[(\xi_h-\overline{\xi},0)] = \int_0^1
A(t,\xi_h,y_h) \xi_h' \cdot (\xi_h'-\overline{\xi}')\, dt +
\frac{1}{2} \int_0^1 \Bigl[ D_2
  A(t,\xi_h,y_h)[\xi_h-\overline{\xi}]\xi_h' \cdot \xi_h' \\
  + B(t,\xi_h,y_h)\cdot (\xi_h' - \overline{\xi}') + D_2 
  B(t,\xi_h,y_h) [\xi_h - \overline{\xi}] \cdot \xi_h' +  D_2 C(t,\xi_h,y_h)
[\xi_h - \overline{\xi}]\\ - D_2 \tilde{\beta} (t,\xi_h,y_h)[\xi_h -
  \overline{\xi}] |y_h'|^2 \Bigr] \, dt,
\end{eqnarray*}
and from the fact that $\xi_h\rightarrow \overline{\xi}$ uniformly and
$\xi_h' \rightarrow \overline{\xi}'$ weakly in $L^2$, we deduce that
the first integral above is infinitesimal. Therefore, also
\[
\int_0^1 A(t,\xi_h,y_h) (\xi_h'-\overline{\xi}') \cdot
(\xi_h'-\overline{\xi}')\, dt
\]
is infinitesimal, and the fact that $A$ is positive definite implies
that $\xi_h \rightarrow \overline{\xi}$ in $W^{1,2}$. Similarly, $D
\tilde{E} (\xi_h,y_h)[(0,y_h - \overline{y})]$ is infinitesimal, so
the last integral in the expression
\begin{eqnarray*}
D\tilde{E}(\xi_h,y_h)[(0,y_h - \overline{y})] = \frac{1}{2} \int_0^1
\Bigl[ D_3 A(t,\xi_h,y_h)[y_h-\overline{y}]\xi_h' \cdot \xi_h' + D_3
  B(t,\xi_h,y_h) [y_h -\overline{y}]\cdot \xi_h' \\
+ D_3 C(t,\xi_h,y_h)[y_h - \overline{y}] - D_3
\tilde{\beta}(t,\xi_h,y_h) [y_h - \overline{y}] |y_h'|^2 \Bigr] \, dt
- \int_0^1 \tilde{\beta} (t,\xi_h,y_h) y_h' (y_h'-\overline{y}')\, dt
\end{eqnarray*} 
is infinitesimal. Therefore, also
\[
\int_0^1 \tilde{\beta} (t,\xi_h,y_h) |y_h'-\overline{y}'|^2\, dt
\]
is infinitesimal, and the fact that $\tilde{\beta}$ is strictly
positive implies that $(y_h - \overline{y})$ converges to $0$ in $W^{1,2}$,
concluding the proof.
\end{proof} \qed 

\section{A pseudo-gradient vector field for the energy functional}

Let $\chi$ be a smooth real function on $\R$ such that $\chi(s)=0$ for
$s\leq 0$, $\chi(s)=1$ for $s\geq 1$, and $\chi'\geq 0$. We choose two
numbers
\begin{equation}
\label{range}
\lambda:= \frac{2}{p_1(s_0,\underline{\beta},\overline{\beta},a,b)}, 
\quad c_0 < -\frac{1}{2} (\overline{\beta}+\lambda)
q_1(s_0,\underline{\beta}, \overline{\beta},a,b,y_1-y_0).
\end{equation}
By Lemma \ref{stima2} with $\mu=0$, 
\[
\mathrm{crit}\, E_{\alpha,\beta} \subset
\{E_{\alpha,\beta+\lambda}>c_0\}.
\]
By the same lemma with $\mu=\lambda$, the manifold $\mathscr{M}$ is
covered by the open sets $\{\nabla E_{\alpha,\beta+\lambda}\neq 0\}$
and $\{E_{\alpha,\beta+\lambda}>c_0\}$. Therefore, the formula
\[
-G = \nabla E_{\alpha,\beta} + \chi(c_0-E_{\alpha,\beta+\lambda}) \frac{
  \|\nabla E_{\alpha,\beta}\|}{\|\nabla E_{\alpha,\beta+\lambda}\|}
\nabla E_{\alpha,\beta+\lambda}
\]
defines a smooth vector field $G$ on $\mathscr{M}$, which coincides
with $-\nabla E_{\alpha,\beta}$ on $\{E_{\alpha,\beta+\lambda}>c_0\}$. 
It is useful to make $G$ bounded, by multiplying it for a suitable
smooth positive function, obtaining the smooth bounded vector field
\[
F:= \frac{G}{\sqrt{1+\|G\|^2}}.
\]

\begin{lem}
\label{lyapunov}
\begin{enumerate}
\item $F$ is a Morse vector field, and $\mathrm{rest} \, F =
  \mathrm{crit} \, E_{\alpha,\beta}$;
\item $E_{\alpha,\beta}$ is a non-degenerate Lyapunov function for $F$
  on $\mathscr{M}$;
\item $E_{\alpha,\beta+\mu}$ is a Lyapunov function for $F$
on $\{E_{\alpha,\beta+\lambda}<c_0-1\}$ (where $F$ has no rest points).
\end{enumerate}
\end{lem}

\begin{proof}
Since the function $(1+\|G\|^2)^{-1/2}$ is positive, it is enough to
prove (i), (ii), and (iii) for the vector field $G$.
The numbers $DE_{\alpha,\beta} [-G]$ and $DE_{\alpha,\beta+\lambda}
[-G]$ are 
\begin{eqnarray}
\label{qua1}
DE_{\alpha,\beta} [-G] = \langle \nabla E_{\alpha,\beta}, -G \rangle =
\|\nabla E_{\alpha,\beta}\|^2 + \chi(c_0-E_{\alpha,\beta+\lambda})
\frac{ \|\nabla E_{\alpha,\beta}\|}{ \|\nabla
  E_{\alpha,\beta+\lambda}\|} \langle \nabla E_{\alpha,\beta} , \nabla
E_{\alpha,\beta+\lambda} \rangle, \\
\label{qua2}
DE_{\alpha,\beta+\lambda} [-G] = \langle \nabla
E_{\alpha,\beta+\lambda}, -G \rangle = \langle \nabla
E_{\alpha,\beta}, \nabla E_{\alpha,\beta+\lambda} \rangle +
\chi(c_0 -E_{\alpha,\beta+\lambda}) \| \nabla E_{\alpha,\beta} \| \, \|
\nabla E_{\alpha,\beta+\lambda} \|.
\end{eqnarray}
Since $0\leq \chi\leq 1$, by the Cauchy-Schwarz inequality the quantity
(\ref{qua1}) is greater or equal than zero, and it equals zero if and
only if $\nabla E_{\alpha,\beta}=0$ or
\begin{equation}
\label{dno}
\chi(c_0 -E_{\alpha,\beta+\lambda})=1 \quad \mbox{and} \quad \langle \nabla
E_{\alpha,\beta}, \nabla E_{\alpha,\beta+\lambda} \rangle = - \|
\nabla E_{\alpha,\beta} \| \, \| \nabla E_{\alpha,\beta+\lambda} \|. 
\end{equation}
The latter cannot occur, for it implies $\nabla
E_{\alpha,\beta+\lambda} = - \nu \nabla E_{\alpha,\beta}$ with
$\nu\geq 0$, that is $\nabla E_{\alpha,\beta+\lambda/(1+\nu)} =0$, so
by Lemma \ref{stima2} and by our choice of $c_0$, 
$E_{\alpha,\beta+\lambda}>c_0$, hence $\chi(c_0-E_{\alpha,\beta+\lambda})=0$.

This shows that $DE_{\alpha,\beta}[G]\leq 0$, and it equals zero
exactly on the critical set of $E_{\alpha,\beta}$. Claims (i) and (ii)
follow from the fact that $G=-\nabla E_{\alpha,\beta}$ on a
neighborhood of such a set, and $E_{\alpha,\beta}$ is a Morse function.

On the set $\{E_{\alpha,\beta+\lambda}<c_0-1\}$ the value of 
$\chi(c_0-E_{\alpha,\beta+\lambda})$ is $1$, so by the Cauchy-Schwarz
inequality the quantity (\ref{qua2}) is greater or equal than zero,
and it cannot vanish because as observed, (\ref{dno}) cannot
occur. Therefore $E_{\alpha,\beta+\lambda}$ is a Lyapunov function for
$G$ on $\{E_{\alpha,\beta+\lambda}<c_0-1\}$.
This proves (iii). 
\end{proof} \qed
 
Denote by $P: T \mathscr{M} \rightarrow H$ the map
\[
(x,y,\xi,\eta) \mapsto \eta, \quad \forall (x,y) \in \mathscr{X}
\times \mathscr{Y}, \; \forall \xi\in T_x \mathscr{X}, \; 
\forall \eta\in T_y \mathscr{Y} = H,
\]
that is the differential of the projection of $\mathscr{X} \times 
\mathscr{Y}$ onto the second factor.

\begin{lem} 
\label{c3}
\begin{enumerate} \item The vector field $F$ is Lipschitz continuous 
on bounded subsets of $\mathscr{M}$. 
\item For every $y\in \mathscr{Y}$, the map 
$\mathscr{X} \rightarrow H$, $x \mapsto P F(x,y)$, is compact.
\end{enumerate}
\end{lem}

\begin{proof}
Let $\mathscr{A}$ be a bounded subset of $\mathscr{M}$.
The vector fields $\nabla E_{\alpha,\beta}$ and $\nabla
E_{\alpha,\beta+\lambda}$ are readily seen to be Lipschitz continuous
on $\mathscr{A}$. The function $\chi(c_0 - E_{\alpha,\beta+\lambda}) 
/ \|\nabla
E_{\alpha,\beta+ \lambda}\|$ is also Lipschitz on $\mathscr{A}$. Indeed, it is
enough to check that $\|\nabla E_{\alpha,\beta+ \lambda}\|$ is bounded
away from zero on $\mathscr{A}\cap \{\chi(c_0 - E_{\alpha,\beta+
\lambda})>0\}$,
and this follows from Lemma \ref{psobs} and from the fact that
$E_{\alpha,\beta+\lambda}$ has no critical points in
$\{E_{\alpha,\beta+\lambda} \leq c_0\}$ (again by Lemma \ref{stima2} with
$\mu=\lambda$). We conclude that $G$ is Lipschitz continuous on $B$,
and so is $F$. This proves (i).

Let $\mathscr{X}_0$ be a bounded subset of $\mathscr{X}$, and let 
$y\in \mathscr{Y}$. 
We must show that $PF(\mathscr{X}_0 \times \{y\})$ has compact
closure in $H$. Since $\chi(c_0-E_{\alpha,\beta+\lambda}) 
\|\nabla E_{\alpha,\beta}\|/ \|\nabla E_{\alpha,\beta+\lambda}\|$ is
bounded on bounded sets, it is enough to show that $P\nabla
E_{\alpha,\beta}( \mathscr{X}_0 \times \{y\})$ (and similarly $P\nabla
E_{\alpha,\beta+\lambda}(\mathscr{X}_0 \times \{y\})$) has compact closure in
$H$. Equivalently, we must show that if $(x_n)\subset \mathscr{X}_0$, and
$(v_n)\subset H$ converges weakly to $0$, then
\[
\langle P\nabla E_{\alpha,\beta}(x_n,y) , v_n \rangle
\rightarrow 0.
\]
Up to a subsequence, we may assume that $x_n$ converges uniformly to
some $x\in \mathscr{X}$, and that $v_n$ converges uniformly to $0$. Now, 
\begin{eqnarray*}
\langle P\nabla E_{\alpha,\beta}(x_n,y) , v_n \rangle = \langle
\nabla E_{\alpha,\beta}(x_n,y) , v_n \rangle =
DE_{\alpha,\beta}(x_n,y)[v_n] 
\\ = \frac{1}{2} \int_0^1 g(\partial_y \alpha(x_n,y)
x_n',x_n')v_n \, dt - \frac{1}{2}
\int_0^1 \partial_y \beta(x_n,y) |y'|^2 v_n \, dt
- \int_0^1 \beta(x_n,y) y' v_n'\, dt.
\end{eqnarray*}
The sequences $g(\partial_y \alpha(x_n,y)x_n',x_n')$ and 
$\partial_y \beta(x_n,y) |y'|^2$ are bounded in
$L^1$, so the first two integrals are infinitesimal. The sequence
$\beta(x_n,y)$ is bounded in $L^{\infty}$, so the third integral is
also infinitesimal. This proves (ii).
\end{proof} \qed   

The proof of the Palais-Smale property for $(F,E_{\alpha,\beta})$ 
on bounded subsets of $\mathscr{M}$ makes use of the above lemma,
together with the following elementary inequality for vectors in a 
Hilbert space.  

\begin{lem}
\label{hilvec}
Let $u$ and $v\neq 0$ be vectors in a Hilbert space, and let $\chi\in
[0,1]$. Then
\[
\|u\|^2 + \chi \frac{\|u\|}{\|v\|} \langle u,v \rangle \geq
\frac{1}{2} \min_{0\leq \theta \leq 1} \| \theta v + (1-\theta) u
\|^2.
\]
\end{lem}

\begin{proof}
If $\langle u,v \rangle \geq 0$ the inequality holds trivially, by
choosing $\theta=0$. So we assume $\langle u, v \rangle<0$. In this
case, $\|u-v\|\neq 0$ and
\[
0 \leq \langle u, u-v \rangle = \|u\|^2 - \langle u,v \rangle \leq
\|u\|^2 - 2 \langle u,v \rangle + \|v\|^2 = \|u-v\|^2,
\]
so the number $\langle u,u-v \rangle /\|u-v\|^2$ is in
$[0,1]$. Therefore,
\begin{equation}
\label{stmin}
\min_{0\leq \theta \leq 1} \|\theta v + (1-\theta) u\|^2 = \min_{0\leq
  \theta \leq 1} \|u - \theta (u-v)\|^2 \leq \left\| u -
  \frac{\langle u,u-v \rangle}{\|u-v\|^2} (u-v) \right\|^2 = \|u\|^2 -
  \frac{ \langle u,u-v \rangle^2}{\|u-v\|^2}.
\end{equation}
Moreover,
\begin{equation}
\label{seg}
\begin{split}
\|u\|^2 - \frac{ \langle u,u-v \rangle^2}{\|u-v\|^2} = \|u\|^2 \left(
1 + \langle \frac{u}{\|u\|} , \frac{u-v}{\|u-v\|} \rangle \right) \left(
1 - \langle \frac{u}{\|u\|} , \frac{u-v}{\|u-v\|} \rangle \right)  \\
\leq 2 \|u\|^2  \left(1 - \langle \frac{u}{\|u\|}, \frac{u-v}{\|u-v\|}
\rangle \right) \leq 2 \|u\|^2 \left( 1 + \langle \frac{u}{\|u\|}, 
\frac{v}{\|v\|} \rangle
\right),
\end{split} \end{equation}
the latter inequality following from
\[
- \langle u,v \rangle \|u-v\| \leq \left| \|u\|^2 - \langle u,v
  \rangle \right| \|v\|,
\]
which is easily verified by taking the squares and using the
Cauchy-Schwarz inequality. By (\ref{stmin}) and (\ref{seg}),
\[
\min_{0\leq \theta \leq 1} \|\theta v + (1-\theta) u\|^2 \leq 2 \left(
\|u\|^2 + \frac{\|u\|}{\|v\|} \langle u,v \rangle \right) \leq 2
\left( \|u\|^2 + \chi \frac{\|u\|}{\|v\|} \langle u,v \rangle \right),
\]
because $\langle u,v \rangle <0$ and $\chi\leq 1$. This concludes the
proof.
\end{proof} \qed

\begin{lem}
\label{ps}
Every bounded sequence $z_n=(x_n,y_n)\in \mathscr{M}$ such that
$DE_{\alpha,\beta}(z_n)[F(z_n)]\rightarrow 0$ has a converging 
subsequence.
\end{lem}

\begin{proof}
Since $G$ is bounded on bounded subsets of $\mathscr{M}$, the
assumption implies that also
\[
- DE_{\alpha,\beta}(x_n,y_n) [ G(x_n,y_n) ] = \left( \|\nabla
  E_{\alpha,\beta} \|^2 + \chi(c_0 - E_{\alpha,\beta+\lambda})
  \frac{\|\nabla E_{\alpha,\beta}\|}{\|\nabla
  E_{\alpha,\beta+\lambda}\|} \langle \nabla E_{\alpha,\beta}, \nabla
  E_{\alpha,\beta+\lambda} \rangle \right) (x_n,y_n)
\]
is infinitesimal, so Lemma \ref{hilvec} implies that there exists
$\theta_n\in [0,1]$ such that
\[
 \| \nabla E_{\alpha,\beta+\theta_n \lambda} \| =
\|\theta_n \nabla E_{\alpha,\beta+\lambda} + (1-\theta_n) \nabla
E_{\alpha,\beta} \| 
\]
is infinitesimal. Up to a subsequence, we may assume that $(\theta_n)$ 
converges to $\theta\in [0,1]$, and
\[
\nabla E_{\alpha,\beta+\theta \lambda} (x_n,y_n) = \nabla
E_{\alpha,\beta+\theta_n \lambda} (x_n,y_n) + (\theta_n -
\theta)\lambda y_n
\]
is also infinitesimal in $W^{1,2}$, because $\|y_n'\|_2$ is
bounded. By Lemma \ref{psobs} the sequence $(x_n,y_n)$ is compact.
\end{proof} \qed

\section{The Morse homology of $\mathbf{E_{\alpha,\beta}}$}

By Lemma \ref{lyapunov} (i) and (ii), $E_{\alpha,\beta}$ is a non-degenerate
Lyapunov function for the smooth Morse vector field $F$ on
$\mathscr{M}$. The vector field $F$ is bounded and Lipschitz on
bounded sets (condition (F1)) by Lemma \ref{c3} (i). By Lemma
\ref{lyapunov} (iii), $E_{\alpha,\beta+\lambda}$ is a Lyapunov
function for $F$ on the set $\mathscr{A} = \{E_{\alpha,\beta+\lambda}
< c_0-1\}$, and $F$ has no rest points on the closure of $\mathscr{A}$
(condition (F2)). The functions $E_{\alpha,\beta}$ and
$E_{\alpha,\beta+\lambda}$ are clearly bounded on every bounded subset
of $\mathscr{M}$ (condition (F3)). By Lemma \ref{interliv} 
applied to $\alpha_0=\alpha_1=\alpha$ and $\beta_0=\beta_1=\beta$, 
the set $\{E_{\alpha,\beta}\leq c\} \cap \{E_{\alpha,\beta+\lambda}
\geq c'\}$ is bounded, for every $c,c'\in \R$ (condition (F4)). 
By Lemma \ref{ps}, every bounded Palais-Smale sequence for
$(F,E_{\alpha,\beta})$ has a converging subsequence (condition
(F5)). Finally, conditions (F6) and (F7) are proved in Lemma \ref{exc1}
and Lemma \ref{c3} (ii), respectively. 
Therefore, the Morse homology $H_*(E_{\alpha,\beta})$ is well defined. 

Let $(\alpha_0,\beta_0)$ and $(\alpha_1,\beta_1)$ be two elements of
$\Gamma(s_0,\underline{\alpha},\overline{\alpha},
\underline{\beta},\overline{\beta},a,b)$ such that
the corresponding energy functionals have only non-degenerate critical
points. If 
\begin{equation}
\label{vicino0}
\max \{ \|\alpha_0-\beta_0\|_{\infty} , \|\alpha_1-\beta_1\|_{\infty}
\} < \frac{\lambda \underline{\alpha}}{\underline{\alpha} + \overline{\beta}},
\end{equation}
Lemma \ref{interliv} implies that the pair $(E_{\alpha_0,\beta_0},
E_{\alpha_1,\beta_1 + \lambda})$ and the pair $(E_{\alpha_1,\beta_1},
E_{\alpha_0,\beta_0 + \lambda})$ satisfy (F8), so the 
homomorphisms
\[
\Phi_{0,1}: H_k(E_{\alpha_0,\beta_0}) \rightarrow H_k(E_{\alpha_1,\beta_1}), \quad
\Phi_{1,0}: H_k(E_{\alpha_1,\beta_1}) \rightarrow H_k(E_{\alpha_0,\beta_0}), 
\]
are well-defined. If
\begin{equation}
\label{vicino1}
\max \{ \|\alpha_0-\beta_0\|_{\infty} , \|\alpha_1-\beta_1\|_{\infty}
\} < \frac{\lambda
\underline{\alpha}}{2\underline{\alpha} + 2\overline{\beta}},
\end{equation}
Lemma \ref{interliv} implies that the set
\[
\{E_{\alpha_0,\beta_0}\leq c\} \cap
\{E_{\alpha_1, \beta_1} + E_{\alpha_1,\beta_1-\lambda} \geq c'\} =
\{E_{\alpha_0,\beta_0}\leq c\} \cap \{E_{\alpha_1,
\beta_1-\lambda/2} \geq c'/2\} 
\]
is bounded. Since $(\alpha_0,\beta_0+\lambda/2)$ and
$(\alpha_1,\beta_1+\lambda)$ belong to
$\Gamma(s_0,\underline{\alpha},\overline{\alpha},\underline{\beta}+\lambda/2,
\overline{\beta}+\lambda,a,b)$, if 
\begin{equation}
\label{vicino2}
\max \{ \|\alpha_0-\beta_0\|_{\infty} , \|\alpha_1-\beta_1\|_{\infty}
\} < \frac{\lambda
  \underline{\alpha}}{2\underline{\alpha} + 2\overline{\beta} + \lambda},
\end{equation}
Lemma \ref{interliv} implies that the set
\[
\{E_{\alpha_1,\beta_1} + E_{\alpha_1,\beta_1+\lambda}\leq c\} \cap
\{E_{\alpha_0, \beta_0+\lambda} \geq c'\} =
\{E_{\alpha_1,\beta_1+\lambda/2}\leq c/2\} \cap \{E_{\alpha_0,
\beta_0+\lambda} \geq c'\} 
\]
is bounded. Therefore, if $(\alpha_0,\beta_0)$ and
$(\alpha_1,\beta_1)$ are so close in the $C^0$ norm that
(\ref{vicino2}) holds - and a fortiori also (\ref{vicino0}) and
(\ref{vicino1}) hold - condition (F9) is satisfied. We conclude that 
the composition
\[
\Phi_{1,0} \circ \Phi_{0,1} : H_k(E_{\alpha_0,\beta_0})
\rightarrow H_k(E_{\alpha_0,\beta_0})
\]
is the identity. Exchanging the role of $(\alpha_0,\beta_0)$ and
$(\alpha_1,\beta_1)$, we deduce that also the composition
\[
\Phi_{0,1} \circ \Phi_{1,0} : H_k(E_{\alpha_1,\beta_1})
\rightarrow H_k(E_{\alpha_1,\beta_1})
\]
is the identity. We conclude that under assumption (\ref{vicino2}), 
the Morse homology of $E_{\alpha_0,\beta_0}$ is isomorphic to the 
Morse homology of $E_{\alpha_1,\beta_1}$. 

The set
$\Gamma(s_0,\underline{\alpha},\overline{\alpha},\underline{\beta},
\overline{\beta},a,b)$ is bounded in $C^0$, and it contains the product
Lorentzian structure $(\underline{\alpha},\underline{\beta})$. 
The set of Lorentzian structures
$(\alpha,\beta)$ in $\Gamma(s_0,\underline{\alpha},\overline{\alpha},
\underline{\beta}, \overline{\beta},a,b)$ such that $E_{\alpha,\beta}$
is a Morse function is $C^0$-dense (actually, $C^{\infty}$-dense). For instance, this follows from the fact that given a globally hyperbolic Lorentzian manifold $(M,h)$ and a point $z_0\in M$, the set of points $z_1\in M$ which are non-conjugate to $z_0$ along every geodesic is residual\footnote{Indeed, by the above fact we can find a sequence of diffeomorphisms $(\varphi_n)$ of $X\times \R$ of the form $(x,y) \mapsto (\varphi_n^1(x),\varphi_n^2(y))$ which converge to the identity in the $C^{\infty}$ topology, and such that $\varphi_n(z_0)$ and $\varphi_n(z_1)$ are non-conjugate along every geodesic. Pulling back the Lorentzian metric by this sequence of diffeomorphisms, we obtain a sequence of globally hyperbolic Lorentzian metrics $C^{\infty}$-converging to the the original one, and
for which $z_0$ and $z_1$ are non-conjugate along every geodesic.}
 (see \cite{uhl75}, Theorem 1 (a)). 
These facts imply that the Morse homology of
every Morse functional $E_{\alpha,\beta}$, $(\alpha,\beta)\in
\Gamma(s_0,\underline{\alpha},\overline{\alpha},\underline{\beta},
\overline{\beta},a,b)$, is isomorphic to the Morse homology of
$E_{\underline{\alpha}, \underline{\beta}}$ (which we may assume to be
a Morse functional on $\mathscr{M}$, by perturbing the metric $g$).   

The Morse homology of $E_{\underline{\alpha},\underline{\beta}}$ is
isomorphic to the singular homology of $\mathscr{X}$ (see Remark
\ref{montreal}). The conclusion follows from the fact that
$\mathscr{X}$ is homotopically equivalent to $\Omega(X)$, the based
loop space of $X$. 

\providecommand{\bysame}{\leavevmode\hbox to3em{\hrulefill}\thinspace}
\providecommand{\MR}{\relax\ifhmode\unskip\space\fi MR }
\providecommand{\MRhref}[2]{%
  \href{http://www.ams.org/mathscinet-getitem?mr=#1}{#2}
}
\providecommand{\href}[2]{#2}

\end{document}